\documentclass{amsart}

\newtheorem{thm}{Theorem}[section]
\newtheorem{cor}[thm]{Corollary}

\newtheorem{prop}[thm]{Proposition}
\theoremstyle{definition}
\newtheorem{defn}[thm]{Definition}
\theoremstyle{remark}
\newtheorem{rem}[thm]{Remark}
\newtheorem{ex}[thm]{Example}

\begin{document}
\title[Some constructions of para-hyperhermitian structures]{\bf Some constructions
of almost para-hyperhermitian structures on manifolds and tangent
bundles}

\author[S. Ianu\c s, G. E. V\^\i lcu]{Stere Ianu\c s, Gabriel Eduard V\^\i lcu}
\thanks{The authors are partially supported by
grant D11-22 CEEX 2006-2008. This paper has been presented in the
"4th German-Romanian Seminar on Geometry" Dortmund, Germany, 15-18
July 2007.}

\date{}
\maketitle \abstract In this paper we give some examples of almost
para-hyperhermitian structures on the tangent bundle of an almost
product manifold, on the product manifold $M\times\mathbb{R}$, where
$M$ is a manifold endowed with a mixed 3-structure and
on the circle bundle over a manifold with a mixed 3-structure.

\endabstract

\section{Introduction}

The para-hyperhermitian structures arise in a natural way in
theoretical physics, both in string theory and integrable systems
(\cite{CMMS}, \cite{DNJ}, \cite{DW}, \cite{H}, \cite{OV}). This kind
of structures have been intensively studied with different names by
many authors (see \cite{AK}, \cite{HIT}, \cite{IZ}, \cite{KMD},
\cite{LBM}, \cite{V} and more).

In this paper we construct some classes of manifolds endowed with
almost para-hyperhermitian structures. The paper is organized as
follows. In section 2 we recall the definition and fundamental
properties of almost para-hyperhermitian manifolds.

In Section 3 we give an almost para-hyperhermitian structure on the
tangent bundle of an almost para-hermitian manifold and study its
integrability.

The concept of mixed 3-structure has been introduced in \cite{IMV}.
In section 4 we study the manifolds endowed with such structures and
prove that $M$ is a manifold with a mixed 3-structure, then
$\overline{M}=M\times\mathbb{R}$ can be endowed with an almost
para-hyperhermitian structure. We construct also an almost
para-hyperhermitian structures on a principal circle bundle over a
smooth manifold endowed with a metric mixed 3-structure.

\section{Preliminaries on almost para-hyperhermitian
manifolds}

An almost product structure on a smooth manifold $M$ is a tensor
field $P$ of type (1,1) on $M$, $P\neq\pm Id$, such that:
\begin{equation}\label{1}
         P^2=Id.
         \end{equation}

An almost complex structure on a smooth manifold $M$ is a tensor
field $J$ of type (1,1) on $M$ such that:
\begin{equation}\label{2}
         J^2=-Id.
         \end{equation}

An almost para-hypercomplex structure on a smooth manifold $M$ is a
triple $H=(J_{\alpha})_{\alpha=\overline{1,3}}$, where $J_1$ is an
almost complex structure on $M$ and $J_2$, $J_3$ are almost product
structures on $M$, satisfying:
\begin{equation}\label{3}
         J_2J_1=-J_1J_2=J_3.
         \end{equation}
In this case $(M,H)$ is said to be an almost para-hypercomplex
manifold.

It is easy to see that any almost para-hypercomplex manifold is of
dimension $4n$, for a non-negative integer $n$.

A semi-Riemannian metric $g$ on $(M,H)$ is said to be
para-hyperhermitian if it satisfies:
 \begin{equation}\label{4}
         g(J_\alpha X,J_\alpha Y)=\epsilon_{\alpha}
         g(X,Y),\forall\alpha=\overline{1,3}
         \end{equation}
for all vector fields $X$,$Y$ on $M$, where $\epsilon_1=1,
\epsilon_2=\epsilon_3=-1$.

Moreover, $(M,g,H)$ is said to be an almost para-hyperhermitian
manifold. It is clear that the signature of $g$ is $(2n,2n)$
because, at each point of $M$, there is a local pseudo-orthonormal
frame field, called adapted frame, given in the following way:
$$\{E_1,...,E_m,J_1E_1,...,J_1E_m,J_2E_1,...,J_2E_m,J_3E_1,...,J_3E_m\}.$$

If $h$ is an arbitrary semi-Riemannian metric on an almost
para-hypercomplex manifold $(M,H)$, then we can define always a
para-hyperhermitian metric $g$ on $M$ by:
\begin{equation}\label{5}
         g(X,Y)=\frac{1}{4}[h(X,Y)+\sum_{\alpha=1}^{3}\epsilon_{\alpha}
         h(J_\alpha X,J_\alpha Y)]
         \end{equation}
for all vector fields $X$,$Y$ on $M$.

An almost para-hypercomplex manifold $(M,H)$ is said to be a
para-hypercomplex manifold if each $J_{\alpha}$, $\alpha=1,2,3$, is
integrable, that is, if the corresponding Nijenhuis tensors:
\begin{equation}\label{6}
         N_{\alpha}(X,Y)=[J_{\alpha}X, J_{\alpha}Y]-J_{\alpha}[X,J_{\alpha}Y]
         -J_{\alpha}[J_{\alpha}X,Y]-\epsilon_{\alpha}[X,Y]
         \end{equation}
$\alpha=1,2,3$, vanish for all vector fields $X$,$Y$ on $M$. In this
case $H$ is said to be a para-hypercomplex structure on $M$.

We remark that if two of the structures ${J_1,J_2,J_3}$ are
integrable, then the third structure is also integrable because we
have:
\begin{eqnarray}\label{7}
          2N_{\alpha}(X,Y)&=&N_{\beta}(J_{\gamma}X,J_{\gamma}Y)+N_{\gamma}(J_{\beta}X,J_{\beta}Y)
          -J_{\beta}N_{\gamma}(J_{\beta}X,Y)-J_{\beta}N_{\gamma}(X,J_{\beta}Y)\nonumber\\
          &&-J_{\gamma}N_{\beta}(J_{\gamma}X,Y)-J_{\gamma}N_{\beta}(X,J_{\gamma}Y)
          +\epsilon_{\alpha}\epsilon_{\beta}N_{\beta}(X,Y)+\epsilon_{\alpha}\epsilon_{\gamma}N_{\gamma}(X,Y)
          \end{eqnarray}
for any even permutation $(\alpha,\beta,\gamma)$ of (1,2,3), where
$\epsilon_1=1, \epsilon_2=\epsilon_3=-1$.

An almost hermitian paracuaternionic manifold is a triple
$(M,\sigma,g)$, where $M$ is a smooth manifold, $\sigma$ is a rank
3-subbundle of $End(TM)$ which is locally spanned by an almost
para-hypercomplex structure $H=(J_{\alpha})_{\alpha=\overline{1,3}}$
and $g$ is a para-hyperhermitian metric in respect with $H$.

\section{An almost para-hyperhermitian structure on the tangent bundle of an almost para-hermitian manifold}

An almost para-hermitian structure on a differentiable manifold $M$
is a pair $(P,g)$, where $P$ is an almost product structure on $M$
and $g$ is a semi-Riemannian metric on $M$ satisfying:
\begin{equation}\label{8}
         g(PX,PY)=-g(X,Y),
         \end{equation}
for all vector fields $X$,$Y$ on $M$.

In this case, $(M,P,g)$ is said to be an almost para-hermitian
manifold. It is easy to see that the dimension of $M$ is even.
Moreover, if $\nabla P=0$ then $(M,P,g)$ is said to be a
para-K\"{a}hler manifold.

We remark that this kind of manifolds appeared for the first time in
\cite{LBM}.

Now, let $(M,P,g)$ be an almost para-hermitian manifold and $TM$ be
the tangent bundle, endowed with the Sasakian metric:
$$G(X,Y)=(g(KX,KY)+g(\pi_*X,\pi_*Y))\circ\pi$$
for all vector fields $X,Y$ on $TM$, where $\pi$ is the natural
projection of $TM$ onto $M$ and $K$ is the connection map (see
\cite{D}).

We remark that if $X\in\Gamma(TM)$, then there exists exactly one
vector field on $TM$ called the "horizontal lift" (resp. "vertical
lift") of $X$ such that for all $t\in TM$:
$$\pi_*X_t^h=X_{\pi(t)},\ \pi_*X_t^v=0_{\pi(t)},\ KX_t^h=0_{\pi(t)},\ KX_t^v=X_{\pi(t)}.$$

\begin{rem}\label{3.1}
It is immediately that for any affine connection $\nabla$ on $M$ we
have:
\begin{eqnarray}
&&[X^h,Y^h]=[X,Y]^h-(R(X,Y)t)^v,\ [X^v,Y^v]=0, \nonumber\\
&&[X^h,Y^v]=(\nabla_XY)^v,\ [X^v,Y^h]=-(\nabla_YX)^v,\nonumber
\end{eqnarray}
for all $X,Y\in T_t(TM)$, $t\in TM$, where $R$ is the curvature
tensor of $\nabla$ on $M$.
\end{rem}

\begin{thm}\label{3.2}
Let $(M,P,g)$ be an almost para-hermitian manifold. Then $TM$ admits
an almost para-hypercomplex structure $H$ which is
para-hyperhermitian in respect to G.
\end{thm}
\begin{proof}
We define three tensor fields $J_1,J_2,J_3$ on $TM$ by the
equalities:
\begin{equation}
J_1 X^h=X^v, J_1 X^v=-X^h,
\end{equation}
\begin{equation}
J_2 X^h=(PX)^v, J_2 X^v=(PX)^h,
\end{equation}
\begin{equation}
J_3 X^h=(PX)^h,J_3 X^v=-(PX)^v.
\end{equation}

We can easily see that we have:
$$J_1^2=-J_2^2=-J_3^2=-Id,$$
$$J_2J_1=-J_1J_2=J_3$$
and
$$G(J_1X,J_1Y)=-G(J_2X,J_2Y)=-G(J_3X,J_3Y)=G(X,Y).$$
\end{proof}

\begin{cor}
The tangent bundle of the pseudosphere $S^6_3$ admits an almost
para-hyperhermitian structure.
\end{cor}
\begin{proof}
The assertion is clear from above Theorem because the pseudosphere
$S^6_3$ can be endowed with an almost para-hermitian structure (see
\cite{BEJ}).
\end{proof}

\begin{cor}
Let $(M,g)$ be a semi-Riemannian manifold and $T^*M$ its cotangent
bundle. Then the tangent bundle of the cotangent bundle $TT^*M$
admits an almost para-hyperhermitian structure.
\end{cor}
\begin{proof}
The assertion is clear from Theorem \ref{3.2} because the cotangent
bundle $T^*M$ can be endowed with an almost para-hermitian structure
(see \cite{OPM}).
\end{proof}

\begin{thm}\label{3.5}
Let $(M,P,g)$ be an almost para-hermitian manifold. Then the almost
para-hypercomplex structure $H=(J_{\alpha})_{\alpha=\overline{1,3}}$
on $TM$ given by Theorem \ref{3.2} is integrable if and only if
$(M,P)$ is a flat para-K\"{a}hler manifold.
\end{thm}
\begin{proof}
Using Remark \ref{3.1} we deduce that the Nijenhuis tensor of $J_1$
is given by:
\begin{eqnarray}\label{12}
N_1(X^h,Y^h)&=&[J_1X^h,J_1Y^h]-J_1[X^h,J_1Y^h]-J_1[J_1X^h,Y^h]-[X^h,Y^h]\nonumber\\
&=&J_1(\nabla_YX)^v-J_1(\nabla_XY)^v-[X,Y]^h+(R(X,Y)t)^v\nonumber\\
&=&(R(X,Y)t)^v,
\end{eqnarray}
\begin{eqnarray}\label{13}
N_1(X^v,Y^v)&=&[J_1X^v,J_1Y^v]-J_1[X^v,J_1Y^v]-J_1[J_1X^v,Y^v]-[X^v,Y^v]\nonumber\\
&=&[X,Y]^h-(R(X,Y)t)^v-J_1(\nabla_YX)^v+J_1(\nabla_XY)^v\nonumber\\
&=&-(R(X,Y)t)^v,
\end{eqnarray}
\begin{eqnarray}\label{14}
N_1(X^h,Y^v)&=&[J_1X^h,J_1Y^v]-J_1[X^h,J_1Y^v]-J_1[J_1X^h,Y^v]-[X^h,Y^v]\nonumber\\
&=&(\nabla_YX)^v+J_1([X,Y]^h-(R(X,Y)t)^v)-(\nabla_XY)^v\nonumber\\
&=&(R(X,Y)t)^h,
\end{eqnarray}
\begin{eqnarray}\label{15}
N_1(X^v,Y^h)&=&[J_1X^v,J_1Y^h]-J_1[X^v,J_1Y^h]-J_1[J_1X^v,Y^h]-[X^v,Y^h]\nonumber\\
&=&-(\nabla_XY)^v+J_1([X,Y]^h-(R(X,Y)t)^v)+(\nabla_YX)^v\nonumber\\
&=&(R(X,Y)t)^h.
\end{eqnarray}
for all $X,Y\in T_t(TM)$, $t\in TM$.

Similarly as above, we deduce:
\begin{equation}\label{16}
N_2(X^h,Y^h)=(P(\nabla_Y P)X-P(\nabla_X P)Y)^h-(R(X,Y)t)^v,
\end{equation}
\begin{equation}\label{17}
N_2(X^v,Y^v)=((\nabla_{PX} P)Y-(\nabla_{PY} P)X)^h-(R(PX,PY)t)^v,
\end{equation}
\begin{equation}\label{18}
N_2(X^h,Y^v)=-(P(\nabla_X P)Y+(\nabla_{PY} P)X)^v+(PR(X,PY)t)^h,
\end{equation}
\begin{equation}\label{19}
N_2(X^v,Y^h)=((\nabla_{PX} P)Y+P(\nabla_Y P)X)^v+(PR(PX,Y)t)^h
\end{equation}
and
\begin{eqnarray}\label{20}
N_3(X^h,Y^h)&=&((\nabla_X P)PY-(\nabla_{PY} P)X+(\nabla_{PX}
P)Y+P(\nabla_{Y}P)X)^h\nonumber\\
             &&-(R(X,Y)t+R(PX,PY)t+PR(PX,Y)t+PR(X,PY)t)^v,
\end{eqnarray}
\begin{equation}\label{21}
N_3(X^v,Y^v)=0,
\end{equation}
\begin{equation}\label{22}
N_3(X^h,Y^v)=(-(\nabla_{PX} P)Y+(\nabla_X P)(PY))^v,
\end{equation}
\begin{equation}\label{23}
N_3(X^v,Y^h)=((\nabla_{PY} P)X-(\nabla_Y P)(PX))^v.
\end{equation}

The proof is now complete from (\ref{12})-(\ref{23}).
\end{proof}

\begin{cor}
Let $(M,P,g)$ be an almost para-hermitian manifold and
$H=(J_{\alpha})_{\alpha=\overline{1,3}}$ the almost
para-hypercomplex structure on $TM$ given in Theorem \ref{3.2}. If
$J_2$ is integrable, then $H$ is a para-hypercomplex structure on
$TM$.
\end{cor}
\begin{proof}
The proof is clear from Theorem \ref{3.5}.
\end{proof}

\section{Almost para-hyperhermitian structures and manifolds endowed with mixed 3-structures}

\begin{defn}
        Let $M$ be a differentiable manifold equipped with a triple
        $(\phi,\xi,\eta)$, where $\phi$ is a a field  of endomorphisms
        of the tangent spaces, $\xi$ is a vector field and $\eta$ is a
        1-form on $M$. If we have:
\begin{equation}\label{24}
        \phi^2=-\epsilon I+\eta\otimes\xi,\  \eta(\xi)=\epsilon
\end{equation}
        then we say that:\\
        i. $(\phi,\xi,\eta)$ is an almost contact
        structure on $M$, if $\epsilon=1$ (cf. \cite{BLR}).\\
        ii.$(\phi,\xi,\eta)$ is a Lorentzian almost paracontact
        structure on $M$, if $\epsilon=-1$ (cf. \cite{MTS}).
\end{defn}

We remark that many authors also include in the above definition the
conditions that $\phi\xi=0$ and $\eta\circ\phi=0$, although these
are deducible from the conditions (\ref{24}) (see \cite{BLR}).

\begin{defn}
        (\cite{IMV})
        Let $M$ be a differentiable manifold which admits an almost contact structure
        $(\phi_1,\xi_1,\eta_1)$ and two Lorentzian almost paracontact
        structures $(\phi_2,\xi_2,\eta_2)$ and
        $(\phi_3,\xi_3,\eta_3)$, satisfying the following
        conditions:\\
        i.
        \begin{equation}\label{25}
        \eta_\alpha(\xi_\beta)=0, \forall \alpha\neq\beta;
        \end{equation}
        ii.
        \begin{equation}\label{26}
        \phi_\alpha(\xi_\beta)=-\phi_\beta(\xi_\alpha)=\epsilon_\gamma\xi_\gamma;
        \end{equation}
        iii.
        \begin{equation}\label{27}
        \eta_\alpha\circ\phi_\beta=-\eta_\beta\circ\phi_\alpha=\epsilon_\gamma\eta_\gamma;
        \end{equation}
        iv.
        \begin{equation}\label{28}
        \phi_\alpha\phi_\beta-\eta_\beta\otimes\xi_\alpha=
        -\phi_\beta\phi_\alpha+\eta_\alpha\otimes\xi_\beta=\epsilon_\gamma\phi_\gamma,
        \end{equation}
        where in (\ref{26}), (\ref{27}) and (\ref{28}), $(\alpha,\beta,\gamma)$ is an even
        permutation of (1,2,3) and $\epsilon_1=1, \epsilon_2=\epsilon_3=-1$.

        Then the manifold $M$ is said to have a mixed 3-structure
$(\phi_\alpha,\xi_\alpha,\eta_\alpha)_{\alpha=\overline{1,3}}$.
\end{defn}

\begin{defn}
If a manifold $M$ with a mixed 3-structure
$(\phi_\alpha,\xi_\alpha,\eta_\alpha)_{\alpha=\overline{1,3}}$
admits a semi-Riemannian metric $g$ such that:
\begin{equation}
g(\phi_\alpha X, \phi_\alpha Y)=\epsilon_\alpha
g(X,Y)-\eta_\alpha(X)\eta_\alpha(Y),
\end{equation}
and
\begin{equation}
g(X,\xi_\alpha)=\eta_\alpha(X)
\end{equation}
for all $X,Y\in\Gamma(TM)$ and $\alpha=1,2,3$, then we say that $M$
has a metric mixed 3-structure and $g$ is called a compatible
metric.

Moreover, if $(\phi_1,\xi_1,\eta_1,g)$ is a Sasakian structure,
i.e.(see \cite{BLR}):
$$(\nabla_X\phi_1) Y=g(X,Y)\xi_1-\eta_1(Y)X$$
and $(\phi_2,\xi_2,\eta_2,g)$, $(\phi_3,\xi_3,\eta_3,g)$ are
Lorentzian para-Sasakian structures, i.e.(see \cite{MTS}):
$$(\nabla_X\phi_2) Y=g(\phi_2X,\phi_2Y)\xi_2+\eta_2(Y)\phi_2^2X,$$
$$(\nabla_X\phi_3) Y=g(\phi_3X,\phi_3Y)\xi_3+\eta_3(Y)\phi_3^2X,$$
where $\nabla$ is the Levi-Civita connection of $g$, then
$((\phi_\alpha,\xi_\alpha,\eta_\alpha)_{\alpha=\overline{1,3}},g)$
is said to be a mixed Sasakian 3-structure on $M$.
\end{defn}
\begin{prop}
Any manifold $M$ with a mixed 3-structure
$(\phi_\alpha,\xi_\alpha,\eta_\alpha)_{\alpha=\overline{1,3}}$
admits a compatible semi-Riemannian metric.
\end{prop}
\begin{proof}
A such type of metric can be constructed from any semi-Riemannian
metric $f$ on $M$, in four step:

i. First we define the metric $u$ by:
$$u(X,Y)=f(\phi_1^2X,\phi_1^2Y)+\eta_1(X)\eta_1(Y)$$
and we can see that $u(X,\xi_1)=\eta_1(X)$.

ii. We define now a new metric $v$ by:
$$v(X,Y)=u(\phi_2^2X,\phi_2^2Y)-\eta_2(X)\eta_2(Y)$$
and we have that $v(X,\xi_\alpha)=\eta_\alpha(X)$, for $\alpha=1,2$.

iii. We define a new metric $h$ by:
$$h(X,Y)=v(\phi_3^2X,\phi_3^2Y)-\eta_3(X)\eta_3(Y)$$
and we can easily see that $h(X,\xi_\alpha)=\eta_\alpha(X)$, for
$\alpha=1,2,3$.

iv. Finally we define the metric $g$ by:
$$g(X,Y)=\frac{1}{4}[h(X,Y)+\sum_{\alpha=1}^{3}\epsilon_\alpha[h(\phi_\alpha X,\phi_\alpha
Y)+\eta_\alpha(X)\eta_\alpha(Y)]$$ and a straightforward
verification shows that $g$ provides a compatible semi-Riemannian
metric on
$(M,(\phi_\alpha,\xi_\alpha,\eta_\alpha)_{\alpha=\overline{1,3}})$.
\end{proof}

\begin{rem}
If
$(M^{4n+3},(\phi_\alpha,\xi_\alpha,\eta_\alpha)_{\alpha=\overline{1,3}},g)$
is a manifold with a metric mixed 3-structure, then it is easy to
see that the signature of $g$ is $(2n+1,2n+2)$ because one can check
that, at each point of $M$, there always exists a pseudo-orthonormal
frame field given in the following way:
$$\{(E_i,\phi_1 E_i, \phi_2
E_i, \phi_3 E_i)_{i=1,n}, \xi_1, \xi_2, \xi_3\}.$$
\end{rem}

\begin{ex}\label{4.5}
1. It is easy to see that if we define
$(\phi_\alpha,\xi_\alpha,\eta_\alpha)_{\alpha=\overline{1,3}}$ in
$\mathbb{R}^3$ by their matrices:
$$\phi_1=\left(%
\begin{array}{ccc}
  0 & 0 & 1 \\
  0 & 0 & 0 \\
  -1 & 0 & 0 \\
\end{array}%
\right),\
\phi_2=\left(%
\begin{array}{ccc}
  0 & 0 & 0 \\
  0 & 0 & 1 \\
  0 & 1 & 0 \\
\end{array}%
\right),\
\phi_3=\left(%
\begin{array}{ccc}
  0 & -1 & 0 \\
  -1 & 0 & 0 \\
  0 & 0 & 0 \\
\end{array}%
\right),
$$
$$
\xi_1=\left(%
\begin{array}{c}
  0 \\
  1 \\
  0 \\
\end{array}%
\right),\
\xi_2=\left(%
\begin{array}{c}
  1 \\
  0 \\
  0 \\
\end{array}%
\right),\
\xi_3=\left(%
\begin{array}{c}
  0 \\
  0 \\
  1 \\
\end{array}%
\right),
$$
$$
\eta_1=\left(%
\begin{array}{ccc}
  0 & 1 & 0 \\
\end{array}%
\right),\
\eta_2=\left(%
\begin{array}{ccc}
  -1 & 0 & 0 \\
\end{array}%
\right),\
\eta_3=\left(%
\begin{array}{ccc}
  0 & 0 & -1 \\
\end{array}%
\right),
$$
then $(\phi_\alpha,\xi_\alpha,\eta_\alpha)_{\alpha=\overline{1,3}}$
is a mixed 3-structure on $\mathbb{R}^3$.

We define now
$(\phi'_\alpha,\xi'_\alpha,\eta'_\alpha)_{\alpha=\overline{1,3}}$ in
$\mathbb{R}^{4n+3}$ by:
$$
\phi'_\alpha=\left(%
\begin{array}{cc}
  \phi_\alpha & 0 \\
  0 & J_\alpha \\
\end{array}%
\right), \xi'_\alpha=\left(%
\begin{array}{c}
  \xi_\alpha \\
  0 \\
\end{array}%
\right)
,\eta'_\alpha=\left(%
\begin{array}{cc}
  \eta_\alpha & 0 \\
\end{array}%
\right),
$$
for $\alpha=1,2,3$, where  $J_1$ is the almost complex structure on
$\mathbb{R}^{4n}$ given by:
\begin{equation}\label{31}
J_1((x_i)_{i=\overline{1,4n}})=(-x_2,x_1,-x_4,x_3,...,-x_{4n-2},x_{4n-3},-x_{4n},x_{4n-1}),
\end{equation}
and $J_2,\ J_3$ are almost product structures on
$\mathbb{R}^{4n}$ defined by:
\begin{equation}\label{32}
J_2((x_i)_{i=\overline{1,4n}})=(-x_{4n-1},x_{4n},-x_{4n-3},x_{4n-2},...,-x_3,x_4,-x_1,x_2),
\end{equation}
\begin{equation}\label{33}
J_3((x_i)_{i=\overline{1,4n}})=(x_{4n},x_{4n-1},x_{4n-2},x_{4n-3},...,x_4,x_3,x_2,x_1).
\end{equation}

Since $J_2J_1=-J_1J_2=J_3$, it is easily checked that
$(\phi'_\alpha,\xi'_\alpha,\eta'_\alpha)_{\alpha=\overline{1,3}}$ is
a mixed 3-structure on $\mathbb{R}^{4n+3}$.

2. Let $(\overline{M},\overline{g})$ be a $(m+2)$-dimensional
semi-Riemannian manifold with index $q\in\{1,2,\dots,m+1\}$ and let
$(M,g)$ be a hypersurface of $\overline{M}$, with
$g=\overline{g}_{|M}$. We say that $M$ is a lightlike hypersurface
of $\overline{M}$ if $g$ is of constant rank $m$ (see \cite{BD2}).

We consider the vector bundle $TM^\perp$ whose fibres are defined
by:
$$T_pM^\perp=\{Y_p\in T_p\overline{M}|\overline{g}_p(X_p,Y_p)=0,\forall X_p\in T_pM\},
\forall p\in M.$$

If $S(TM)$ is the complementary distribution of $TM^\perp$ in $TM$,
which is called the screen distribution, then there exists a unique
vector bundle $ltr(TM)$ of rank 1 over $M$ so that for any non-zero
section $\xi$ of $TM^\perp$ on a coordinate neighborhood $U \subset
M$, there exists a unique section $N$ of $ltr(TM)$ on $U$
satisfying:
$$\overline{g}(N,\xi)=1,\ \overline{g}(N,N)=\overline{g}(W,W)=0,\ \forall W \in
\Gamma(S(TM)_{|U})$$ (see \cite{BD2}).

In a lightlike hypersurface $M$ of an almost hermitian
paraquaternionic manifold $(\overline{M},\overline{g},\sigma)$ such
that $\xi$ and $N$ are globally defined on $M$, there is a mixed
3-structure (see \cite{IMV}).

3. The unit sphere $S^{4n+3}_{2n+1}$ is the canonical example of
manifold with a mixed Sasakian 3-structure. This structure is
obtained by taking $S^{4n+3}_{2n+1}$ as hypersurface of
$(\mathbb{R}^{4n+4}_{2n+2},\overline{g})$. It is easy to see that on
the tangent spaces $T_pS^{4n+3}_{2n+1}$, $p\in S^{4n+3}_{2n+1}$, the
induced metric $g$ is of signature $(2n+1,2n+2)$.

If $(J_\alpha)_{\alpha=\overline{1,3}}$ is the canonical
paraquaternionic structure on the $\mathbb{R}^{4n+4}_{2n+2}$ and $N$
is the unit normal vector field to the sphere, we can define three
vector fields on $S^{4n+3}_{2n+1}$ by:
$$\xi_\alpha=-J_\alpha N,\ \alpha=1,2,3.$$

If $X$ is a tangent vector to the sphere then $J_\alpha X$ uniquely
decomposes onto the part tangent to the sphere and the part parallel
to $N$. Denote this decomposition by:
$$J_\alpha X=\phi_\alpha X+\eta_\alpha(X)N.$$
This defines the 1-forms $\eta_\alpha$ and the tensor fields
$\phi_\alpha$ on $S^{4n+3}_{2n+1}$, where $\alpha=1,2,3$.

Now we can easily see that that
$(\phi_\alpha,\xi_\alpha,\eta_\alpha)_{\alpha=\overline{1,3}}$ is a
mixed Sasakian 3-structure on $S^{4n+3}_{2n+1}$.

4. Since we can recognize the unit sphere $S^{4n+3}_{2n+1}$ as the
projective space $P^{4n+3}_{2n+1}(\mathbb{R})$, by identifying
antipodal points, we have also that $P^{4n+3}_{2n+1}(\mathbb{R})$
admits a mixed Sasakian 3-structure.
\end{ex}

\begin{thm}\label{4.6}
If $M$ is a manifold with a mixed 3-structure
$(\phi_\alpha,\xi_\alpha,\eta_\alpha)_{\alpha=\overline{1,3}}$ then
$\overline{M}=M\times\mathbb{I}$ can be endowed with an almost
para-hyperhermitian structure, where $\mathbb{I}$ is $\mathbb{R}$ or
some open interval in $\mathbb{R}$.
\end{thm}
\begin{proof}
We define three tensor fields $J_1,J_2,J_3$ on $\overline{M}$ by the
equalities:
\begin{equation}\label{34}
J_{\alpha}=\left(%
\begin{array}{cc}
  \phi_\alpha & \frac{\xi_\alpha}{f} \\
  -f\eta_\alpha & 0 \\
\end{array}%
\right), \alpha=1,2,3.
\end{equation}
where $f$ is a positive function on $\mathbb{I}$.

By a straightforward computation one can check:
$$J_1^2=-J_2^2=-J_3^2=-Id,$$
$$J_2J_1=-J_1J_2=J_3$$
and the proof is complete because we can construct now a
para-hyperhermitian metric on $\overline{M}=M\times\mathbb{R}$ from
any arbitrary semi-Riemannian metric.
\end{proof}

\begin{cor}
$S^{4n+3}_{2n+1}\times \mathbb{I}$ and
$P^{4n+3}_{2n+1}(\mathbb{R})\times \mathbb{I}$ can be endowed with
almost para-hyperhermitian structures, where $\mathbb{I}$ is
$\mathbb{R}$ or some open interval in $\mathbb{R}$.
\end{cor}
\begin{proof}
From Example \ref{4.5} we have that $S^{4n+3}_{2n+1}$ and
$P^{4n+3}_{2n+1}(\mathbb{R})$ are mixed 3-Sasakian manifolds and,
consequently, each structure $J_\alpha$ given by (\ref{34}) is
integrable.
\end{proof}

\begin{cor}
Let $M$ be a manifold endowed with a mixed 3-structure
$(\phi_\alpha,\xi_\alpha,\eta_\alpha)_{\alpha=\overline{1,3}}$. Then
the dimension of $M$ is $4n+3$, where $n$ is a non-negative integer.
\end{cor}
\begin{proof}
The assertion is clear from Theorem \ref{4.6}.
\end{proof}

\begin{rem}

If $M^{4n+3}$ is  a manifold endowed with a mixed 3-Sasakian
structure
$((\phi_\alpha,\xi_\alpha,\eta_\alpha)_{\alpha=\overline{1,3}},g)$,
then we can define a para-hyper-K\"{a}hler structure
$\{J_\alpha\}_{\alpha=\overline{1,3}}$ on the cone
$(C(M),\overline{g})=(M\times\mathbb{R}_+,dr^2+r^2g)$, by:
\begin{equation}
\left\{\begin{array}{rcl}
       J_\alpha X&=&\phi_\alpha X-\eta_\alpha(X)\Phi\\
       J_\alpha\Phi&=&\xi_\alpha
       \end{array}\right.
\end{equation}
for any vector field $X\in\Gamma(TM)$ and $\alpha=1,2,3$, where
$\Phi=r\partial_r$ is the Euler field on $C(M)$.

Moreover, conversely, if a cone
$(C(M),\overline{g})=(M\times\mathbb{R}_+,dr^2+r^2g)$ admits a
para-hyper-K\"{a}hler structure
$\{J_\alpha\}_{\alpha=\overline{1,3}}$, then we can identify $M$
with $M\times\{1\}$ and we have a mixed 3-Sasakian structure
$((\phi_\alpha,\xi_\alpha,\eta_\alpha)_{\alpha=\overline{1,3}},g)$
on $M$ given by:
\begin{equation}
\xi_\alpha=J_\alpha(\partial_r),\ \phi_\alpha X=\nabla_X
\xi_\alpha,\ \eta_\alpha(X)=g(\xi_\alpha,X),
\end{equation}
for any vector field $X\in\Gamma(TM)$ and $\alpha=1,2,3$.

Finally, since a para-hyper-K\"{a}hler manifold is Ricci-flat, we
conclude that $M$ is an Einstein space with Einstein constant
$\lambda=4n+2$ (see \cite{BG}).

\end{rem}

\begin{rem}

Let $P=P(M,\pi,S^1)$ be a principal circle bundle over a smooth
manifold $M^{4n+3}$ with a metric mixed 3-structure
$((\phi_\alpha,\xi_\alpha,\eta_\alpha)_{\alpha=\overline{1,3}},g)$.

Let $g'$ be a semi-Riemannian metric on $P$ such that
$\pi:(P,g')\rightarrow (M,g)$ is a semi-Riemannian submersion.
Putting $\mathcal{V}_x=Ker\ \pi_{*x}$, for any $x\in M$, we obtain
an integrable distribution $\mathcal{V}$, which is called vertical
distribution and corresponds to the foliation of $M$ determined by
the fibres of $\pi$. The complementary distribution $\mathcal{H}$ of
$\mathcal{V}$, determined by the semi-Riemannian metric $g'$, is
called horizontal distribution. We have now the decomposition:
$$TP=\mathcal{H}\oplus \mathcal{V}.$$

We recall that the sections of $\mathcal{V}$, respectively
$\mathcal{H}$, are called the vertical vector fields, respectively
horizontal vector fields. An horizontal vector field $X^h$ on $P$ is
said to be basic if $X^h$ is $\pi$-related to a vector field $X'$ on
$M$. It is clearly that every vector field $X$ on $M$ has a unique
horizontal lift $X^h$ to $P$ and $X^h$ is basic.
\end{rem}

\begin{thm}\label{5.1}
There is an almost para-hyperhermitian structure on any principal
circle bundle $P=P(M,\pi,S^1)$  over a manifold $M$ with a metric
mixed 3-structure
$((\phi_\alpha,\xi_\alpha,\eta_\alpha)_{\alpha=\overline{1,3}},g)$.
\end{thm}
\begin{proof}
Making use of the metric mixed 3-structure on $M$ we can define
three tensor fields $J_1,J_2,J_3$ on $P$ by the equalities:
\begin{equation}
\left\{\begin{array}{rcl}
       J_\alpha X^h&=&(\phi_\alpha X)^h+\eta_\alpha(X)\Theta\\
       J_\alpha\Theta&=&-(\xi_\alpha)^h
       \end{array}\right.
\end{equation}
for $\alpha=1,2,3$, where $\Theta$ is the nowhere vanishing vertical
vector field on $P$ which generates the $S^1$ action on $P$.

Now, we can easily see that we have:
$$J_1^2=-J_2^2=-J_3^2=-Id,$$
$$J_2J_1=-J_1J_2=J_3$$
and the proof is now complete since we can construct a
para-hyperhermitian metric from any arbitrary semi-Riemannian metric
on $P$.
\end{proof}

\begin{cor}
$\mathbb{R}^{4n+3}_{2n+1}\times S^1$, $S^{4n+3}_{2n+1}\times S^1$
and $P^{4n+3}_{2n+1}(\mathbb{R})\times S^1$ can be endowed with
almost para-hyperhermitian structures.
\end{cor}
\begin{proof}
The assertion is clear from Example \ref{4.5} and Theorem \ref{5.1}.
\end{proof}

\begin{center}

Stere Ianu\c s \\
{\em    University of Bucharest,\\
        Department of Mathematics,\\
        C.P. 10-119, Post. Of. 10, Bucharest 72200, Rom$\hat a$nia}\\
        e-mail: ianus@gta.math.unibuc.ro
\end{center}

\begin{center}
Gabriel Eduard V\^\i lcu \\
{\em  ''Petroleum-Gas'' University of Ploie\c sti,\\
         Department of Mathematics and Computer Science,\\
         Bulevardul Bucure\c sti, Nr. 39, Ploie\c sti, Rom$\hat a$nia}\\
         e-mail: gvilcu@mail.upg-ploiesti.ro\\
\end{center}

\end{document}